# Transient longitudinal waves in 2D square lattices with Voigt elements under concentrated loading

Nadezhda I. Aleksandrova (ORCID iD: 0000-0002-9540-1514)

*Chinakal Institute of Mining of the Siberian Branch of the Russian Academy of Sciences,*

*Krasnyj pr. 54, Novosibirsk 630091 Russia* (nialex@misd.ru)

**Abstract:** The aim of this article is to study the attenuation of transient low-frequency waves in 2D lattices of point masses connected by Voigt elements, under an antiplane concentrated loading. The emphasis is on obtaining analytical estimates for solutions using methods of asymptotic inversion of the Laplace and Fourier transforms in the vicinity of the quasi-front of infinitely long waves. In addition, the problems under study are solved by a finite difference method. The main result of the article is the asymptotic estimates of low-frequency and high-frequency perturbations in the 2D lattice for long periods of time under a transient load. It is shown that the obtained asymptotic estimates qualitatively and quantitatively agree with the results of numerical calculations.

**Keywords:** 2D lattice, transient wave, analytical solution, numerical simulation, antiplane problem, pulse load

**Highlights**

- Asymptotic estimates are found for perturbations caused by a concentrated pulse load
- The degree of attenuation of perturbations with increasing time is studied analytically
- It is shown that the asymptotic and finite-difference solutions agree well

## 1. Introduction

In [1] it was proposed to take into account the block structure of rocks in mathematical models, namely, to consider a rock mass as a system of nested blocks interconnected by interlayers consisting of weaker, fractured rocks. In [2, 3] it is shown that deformations of a block rock mass, both in statics and in dynamics, occur mainly due to compliance of the interlayers, and it is also noted that the block structure gives rise to pendulum waves, i.e., of the waves with low propagation velocity, low frequency, and low attenuation. Such waves are actually observed in rocks. Many researchers studied various aspects of their propagation, see, e.g., [4 – 8]. Nevertheless, some natural questions that arise in the theory of block media remain unexplored. In this article, we will try to answer some of them.

A periodic lattice of masses connected by springs and dampers is the simplest model of a block rock medium. However, even such a simple model makes it possible to describe the low-frequency waves that arise under impact. A comparison of the results of numerical calculations using this model and experimental data is given, e.g., in [9] in the 1D case and in [10] in the 3D case; and comparison of the results of numerical calculations with analytical solutions is given, e.g., in [11] in the 1D case and in [12] in the 2D case.

In this article, a block medium with viscoelastic interlayers is simulated as a 2D square lattice of point masses connected by Voigt elements in axial and diagonal directions. For antiplane problems of wave propagation in viscoelastic lattices under the action of concentrated transient loads, we find asymptotic solutions and compare them with finite-difference solutions.

## 2. Statement of the problem and its solution

The block medium is modeled by a 2D lattice of point masses connected in parallel by springs and dampers in axial and diagonal directions (Fig. 1). As shown in [12], the equations describing plane motions of the 2D square lattice are equivalent to two linearly independent wave equations, each of which contains only one unknown function. By analogy with the theory of elasticity, one of those equations describes the propagation of transverse waves in the lattice, and the other describes the propagation of longitudinal waves. In this article, we study only the propagation of longitudinal waves in the 2D square lattice.

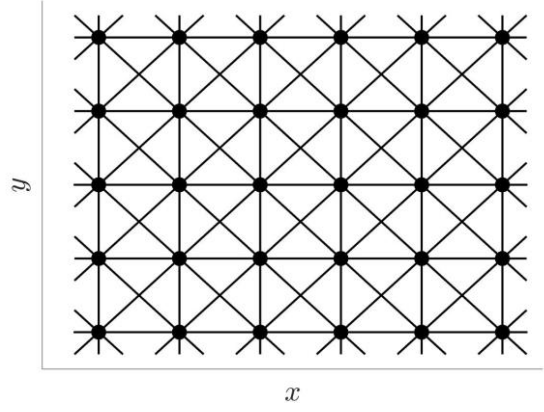

**Fig. 1.** The square lattice of the masses connected by springs in the directions of axes *x*, *y*, and in the diagonal directions.

The motion of the lattice masses is described by the following scalar wave equation with respect to one unknown function $\varphi_{n,m}$:

$$\ddot{\varphi}_{n,m} = D\varphi_{n,m} + \lambda D\dot{\varphi}_{n,m} + Q(t)\delta_{0n}\delta_{0m}. \qquad (1)$$

Here *n*, *m* are the indices of the masses along axes *x*, *y*; $\varphi_{n,m}$ is the displacement of the lattice mass with indices *n*, *m* in the direction orthogonal to the lattice plane; $\lambda$ is the viscosity of a dampers; $\delta_{0n}$ is the Kronecker delta; and $D$ denotes the following operator:

$$D\varphi_{n,m} = \left(\varphi_{n+1,m+1} + \varphi_{n-1,m-1} + \varphi_{n+1,m-1} + \varphi_{n-1,m+1} + \varphi_{n+1,m} + \varphi_{n-1,m} + \varphi_{n,m-1} + \varphi_{n,m+1} - 8\varphi_{n,m}\right)/2.$$

It is also assumed in (1) that *t* is time and $Q(t)$ is the amplitude of a concentrated load, which acts orthogonally to the lattice plane at the point with coordinates $n = 0$, $m = 0$. In (1), the values of point masses, as well as the length and stiffness of the springs, are taken as unity. The initial conditions for (1) are supposed to be zero.

In order to construct an analytical solution, we apply to (1) the Laplace transform with respect to time and the discrete Fourier transforms with respect to coordinates *n*, *m*:

$$f^L(p) = \int_0^\infty f(t)e^{-pt}dt, \quad f(t) = \frac{1}{2\pi i}\int_{\alpha-i\infty}^{\alpha+i\infty} f^L(p)e^{pt}dp,$$

$$f^{F_n}(q_x) = \sum_{n=-\infty}^{n=\infty} f_n e^{iq_x n}, \quad f_n(t) = \frac{1}{2\pi} \int_{-\pi}^{\pi} f^{F_n}(q_x) e^{-iq_x n} dq_x,$$

$$f^{F_m}(q_y) = \sum_{m=-\infty}^{m=\infty} f_m e^{iq_y m}, \quad f_m(t) = \frac{1}{2\pi} \int_{-\pi}^{\pi} f^{F_m}(q_y) e^{-iq_y m} dq_y.$$

Here the superscript $L$ corresponds to the Laplace transform in time with parameter $p$, while the superscripts $F_n$, $F_m$ correspond to the discrete Fourier transforms with respect to coordinates $n$, $m$ with parameters $q_x$, $q_y$, respectively.

The Laplace–Fourier transform of $\varphi$ is as follows:

$$\varphi^{LF_n F_m}(p, q_x, q_y) = \frac{Q^L(p)}{D_1(p, q_x, q_y)}, \qquad (2)$$

where

$$D_1(p, q_x, q_y) = B_1[1 - (2\cos q_x + 1)\cos q_y (\lambda p + 1)/B_1], \quad B_1 = p^2 + (4 - \cos q_x)(\lambda p + 1).$$

Apply to $\varphi^{LF_n F_m}$ the inverse discrete Fourier transform with respect to $q_y$. Using [13], in the case $m = 0$ we get:

$$\varphi_0^{LF_n}(p, q_x) = \frac{Q^L(p)}{\sqrt{6[p^2 + 6\sin^2(q_x/2)(\lambda p + 1)][p^2 + (6 - 2\sin^2(q_x/2))(\lambda p + 1)]}}. \qquad (3)$$

To invert the Laplace and Fourier transforms, we apply different methods in the cases $\lambda = 0$ and $\lambda > 0$.

### 2.1. Step load

Let the dependence of the load on time be described by the Heaviside step function $H(t)$: $Q(t) = H(t)$. In the case $\lambda = 0$, the asymptotic solution of equation (1) with such a load is published in [12]. In particular, the following asymptotic formulas are obtained there for $t \to \infty$ for displacements $\varphi_{n,0}$, velocities $\dot{\varphi}_{n,0}$, and accelerations $\ddot{\varphi}_{n,0}$ of the masses of the lattice in the vicinity of the front of the longitudinal wave:

$$\varphi_{n,0}(t) \sim \frac{1}{2\pi c_1^2} \begin{cases} \ln\left(c_1 t/n + \sqrt{c_1^2 t^2/n^2 - 1}\right) H(c_1 t - n), & n \neq 0, \\ \ln\left(4\sqrt{6} c_1 t\right) + \gamma, & n = 0, \end{cases} \qquad (4)$$

$$\dot{\varphi}_{n,0}(t) \sim J_n^2(c_1 t)/(2c_1) \quad \text{as} \quad n \to \infty, \qquad (5)$$

$$\ddot{\varphi}_{n,0}(t) \sim J_n(c_1 t) J_n'(c_1 t) \quad \text{as} \quad n \to \infty, \qquad (6)$$

where $J_n$ is the Bessel function of the first kind of integer order $n$, the prime denotes the derivative of the function, $\gamma = 0.577...$ is the Euler constant, and $c_1$ is the velocity of infinitely long longitudinal waves in the lattice in the case $\lambda = 0$. Note that in [14] it is shown that the wave velocity $c_1$ is determined by the formula: $c_1 = \sqrt{3/2}$.

We need the following asymptotic formulas for the Bessel function and its derivative, which are valid in a neighborhood of the point $n = c_1 t$:

$$J_n(c_1 t) \sim \frac{\mathrm{Ai}(\kappa)}{(c_1 t/2)^{1/3}} \quad \text{as} \quad t \to \infty \quad \text{and} \quad n \to \infty, \qquad (7)$$

$$J_n'(c_1 t) \sim -\frac{\mathrm{Ai}'(\kappa)}{(c_1 t/2)^{2/3}} \quad \text{as} \quad t \to \infty \quad \text{and} \quad n \to \infty. \qquad (8)$$

Here $\mathrm{Ai}(\kappa)$ is the Airy function and

$$\kappa = \frac{n - c_1 t}{(c_1 t/2)^{1/3}}. \qquad (9)$$

Formula (7) is given in [15] and formula (8) is obtained in [16].

Substituting (7) and (8) into (5) and (6), we obtain

$$\dot{\varphi}_{n,0}(t) \sim \frac{1}{2c_1} \left[\frac{\mathrm{Ai}(\kappa)}{(c_1 t/2)^{1/3}}\right]^2 \quad \text{as} \quad t \to \infty \quad \text{and} \quad n \to \infty, \qquad (10)$$

$$\ddot{\varphi}_{n,0}(t) \sim -\frac{2\mathrm{Ai}(\kappa)\mathrm{Ai}'(\kappa)}{c_1 t} \quad \text{as} \quad t \to \infty \quad \text{and} \quad n \to \infty. \qquad (11)$$

The asymptotic solutions (10), (11) obtained in terms of the Airy function are less accurate than the asymptotic solutions (5), (6) obtained in terms of the Bessel functions. However, the advantage of solutions (10), (11) compared to (5), (6) is that they explicitly describe the degree of attenuation of perturbations with increasing time or distance in the vicinity of the quasi-front $n = c_1 t$ of the longitudinal wave, and describe the degree of expansion of the quasi-front zone, i.e., the zone where perturbations vary from zero to maximum.

In order to determine the limits of applicability of asymptotic solutions (4) – (6), we solve (1) by the finite difference method using an explicit scheme. To do this, the second time derivative in (1) is approximated by the central difference and the first time derivative is approximated by the backward difference. The step of the difference grid in time is denoted by $\tau$. Everywhere below we assume $\tau = 0.01$.

In Figs. 2a – 2c, numerical solutions and asymptotic solutions (4) – (6) are shown for the step load in the case $\lambda = 0$ at the point with coordinates $n = m = 25$. For this, asymptotic solutions (4) – (6) are rewritten in terms of the radial coordinate $r = \sqrt{n^2 + m^2}$ (the corresponding formulas are given explicitly in [12]). In Fig. 2, the vertical dashed lines correspond to the arrival time of the quasi-front of the longitudinal wave, $t = r/c_1$.

Figs. 2b and 2c and Eqs. (9) – (11) show that in the case $\lambda = 0$ for the step load the following effects are observed in the lattice as $t \to \infty$ ($n \to \infty$): the maximum amplitude of the velocities of the masses of the lattice decreases as $t^{-2/3}$ ($n^{-2/3}$); the accelerations of the masses of the lattice decreases as $t^{-1}$ ($n^{-1}$); and the quasi-front zone expands as $t^{1/3}$ ($n^{1/3}$). Besides, Fig. 2a and Eq. (4) show that in this case the amplitudes of the displacements are proportional to $\ln(t/n)$ as $t \to \infty$, $n \neq 0$.

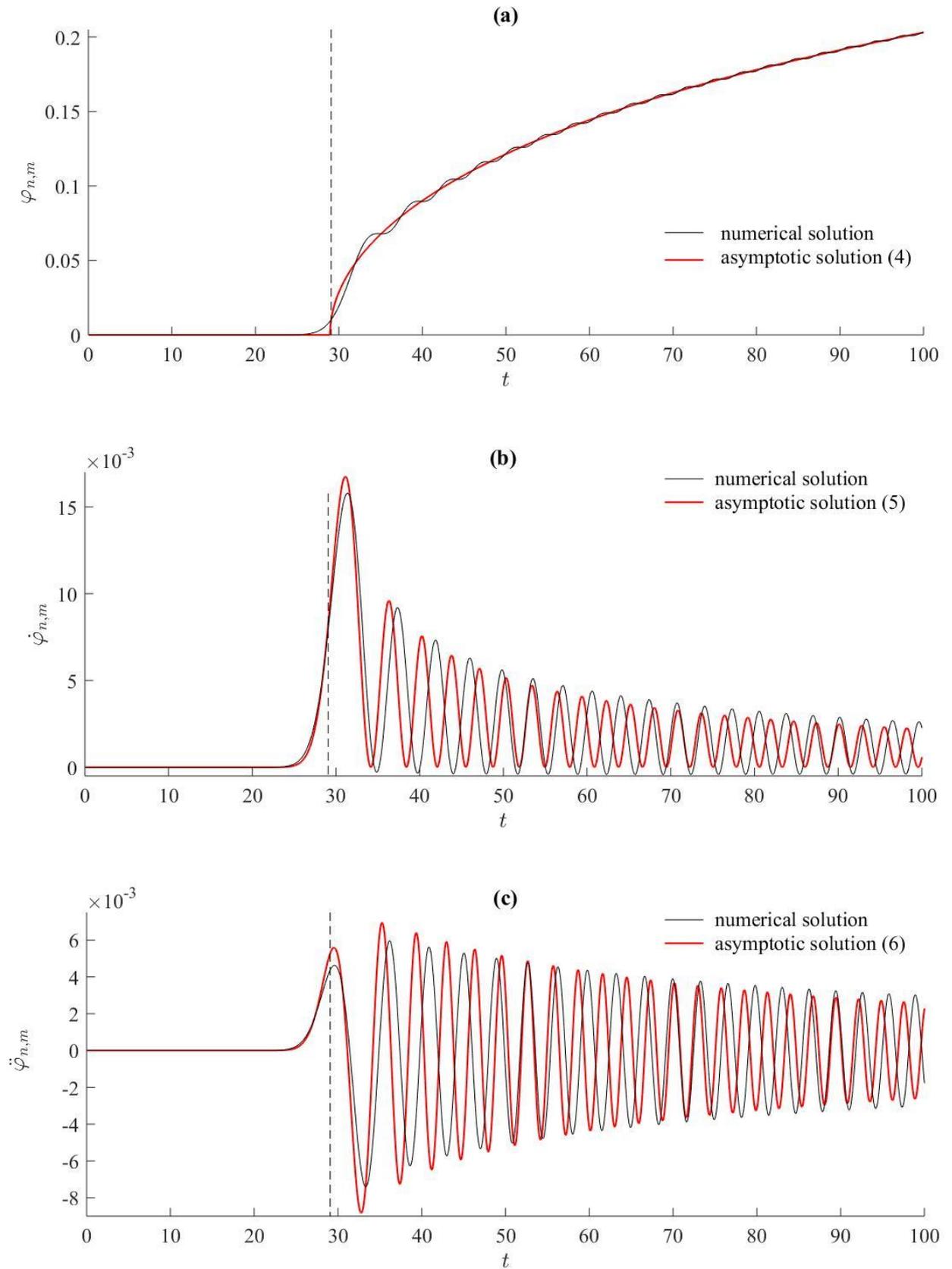

**Fig. 2.** Time dependences of displacements $\varphi_{n,m}$, velocities $\dot{\varphi}_{n,m}$, and accelerations $\ddot{\varphi}_{n,m}$ of the mass with coordinates $n = m = 25$ in the case $\lambda = 0$ for the step load.

Let us continue studying the propagation of waves in a lattice generated by a step load. Now let us consider the case $\lambda > 0$. To find the inverse Fourier transform with respect to $q_x$ and the inverse Laplace transform with respect to $p$ of the function $\dot\varphi_0^{LF_n}$ in (3), we use Slepyan's asymptotic method [17] of joint inversion of these transforms in a vicinity of the ray $n = c_1 t$. To do this, we make the substitution $p = s + iq_x(c_1 + c')$, where $c' = (n - c_1 t)/t \to 0$ and $c'$ defines the vicinity of the ray $n = c_1 t$. We look for the asymptotics of the long-wave perturbations ($|q_x| < \varepsilon$, $\varepsilon$ is small) for large values of time $t$. Recall that the condition $t \to \infty$ in the time domain corresponds to the condition $s \to 0$ in the $s$-domain. Expanding the denominator in (3) into series in powers of the small parameter $q_x$ as $s \to 0$ and $c' \to 0$, and using Slepyan's method [17], we obtain:

$$\dot\varphi_{n,0}(t) = \frac{1}{4\pi^2 i} \int_{-\pi}^{\pi}\int_{\alpha-i\infty}^{\alpha+i\infty} \dot\varphi_0^{LF_n}(s+iq_x(c+c'), q_x) e^{st} ds\, dq_x = \frac{1}{4\pi^2 i} \int_{-\varepsilon}^{\varepsilon}\int_{\alpha-i\infty}^{\alpha+i\infty} \frac{e^{st} ds\, dq_x}{\sqrt{12 iq_x c_1(s + iq_x c' + \alpha q_x^2/t)}},$$

$$\ddot\varphi_{n,0}(t) = \frac{1}{4\pi^2 i} \int_{-\pi}^{\pi}\int_{\alpha-i\infty}^{\alpha+i\infty} \ddot\varphi_0^{LF_n}(s+iq_x(c+c'), q_x) e^{st} ds\, dq_x = \frac{1}{4\pi^2 i} \int_{-\varepsilon}^{\varepsilon}\int_{\alpha-i\infty}^{\alpha+i\infty} \frac{\sqrt{iq_x c_1}\, e^{st} ds\, dq_x}{\sqrt{12(s + iq_x c' + \alpha q_x^2/t)}},$$

where $\alpha = 3\lambda t/4$.

Using [18] and [19], we calculate these integrals. As a result, we have:

$$\dot\varphi_{n,0}(t) = \frac{1}{3\pi^{3/2}(2\lambda t^3)^{1/4}} \Phi_1(\kappa), \tag{12}$$

$$\ddot\varphi_{n,0}(t) = \frac{1}{3\pi^{3/2}(\lambda^3 t^5/2)^{1/4}} \Phi_2(\kappa), \tag{13}$$

where

$$\Phi_1(\kappa) = \int_0^\infty e^{-z^2} z^{-1/2} \sin(-z\kappa + \pi/4)\, dz = 2^{-2}\pi e^{-\eta} |\kappa|^{1/2}\left[I_{-1/4}(\eta) - \operatorname{sgn}(\kappa) I_{1/4}(\eta)\right], \tag{14}$$

$$\Phi_2(\kappa) = \int_0^\infty e^{-z^2} z^{1/2} \sin(z\kappa + \pi/4) dz = \qquad (15)$$
$$= 2^{-4} \pi e^{-\eta} |\kappa|^{3/2} \left\{ I_{-3/4}(\eta) - I_{1/4}(\eta) - \mathrm{sgn}(\kappa) \left[ I_{3/4}(\eta) - I_{-1/4}(\eta) \right] \right\},$$

$$\kappa = \frac{n - c_1 t}{(3\lambda t / 4)^{1/2}}, \qquad (16)$$

and $\eta = \kappa^2/8$. Here, $I_\nu$ is the modified Bessel function [13].

In Figs. 3a – 3c, numerical and asymptotic solutions (4), (12), and (13) are shown for the step load in the case $\lambda = 0.1$ at the point with coordinates $n = m = 25$. For this, asymptotic solutions (4), (12), and (13) are rewritten in terms of the radial coordinate $r = \sqrt{n^2 + m^2}$. In Fig. 3, the vertical dashed lines correspond to the arrival time of the quasi-front of the longitudinal wave, $t = r/c_1$.

Fig. 3b and Fig. 3c and formulas (12), (13), and (16) show that, in the case of the step load and $\lambda > 0$, the following effects are observed in the lattice as $t \to \infty$ ($n \to \infty$): the maximum amplitude of the velocities of the masses of the lattice decreases as $t^{-3/4}$ ($n^{-3/4}$); the maximum amplitude of the accelerations of the masses of the lattice decreases as $t^{-5/4}$ ($n^{-5/4}$); and the quasi-front zone expands as $t^{1/2}$ ($n^{1/2}$). In addition, Eqs. (12), (13), and (16) show that the maximum amplitude of the velocities of the masses of the lattice decreases as $\lambda^{-1/4}$ with an increase in the viscosity parameter $\lambda$; the accelerations of the masses of the lattice decreases as $\lambda^{-3/4}$; and the quasi-front zone expands as $\lambda^{1/2}$.

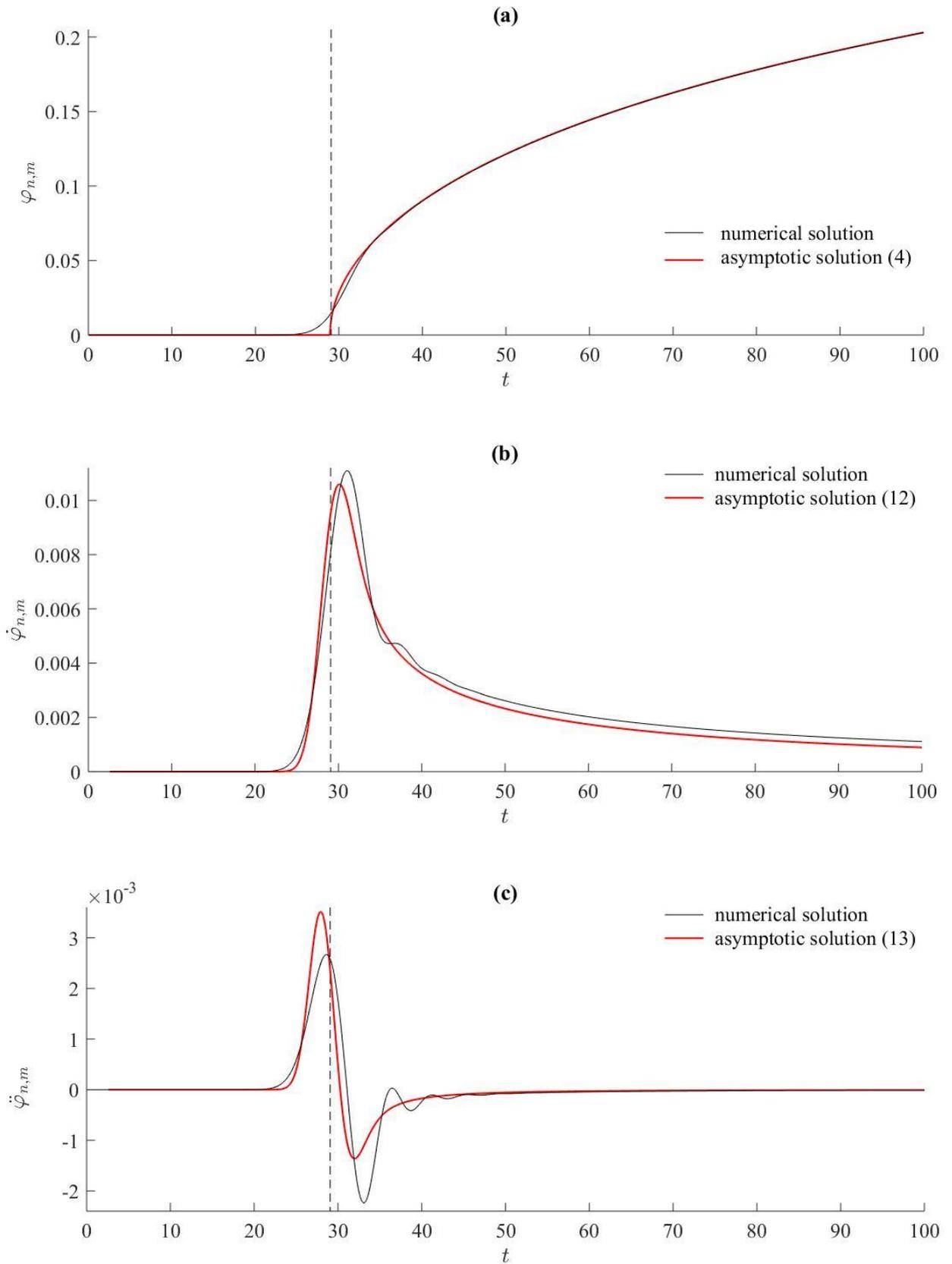

**Fig. 3.** Time dependences of displacements $\varphi_{n,m}$, velocities $\dot{\varphi}_{n,m}$, and accelerations $\ddot{\varphi}_{n,m}$ of the mass with coordinates $n = m = 25$ in the case $\lambda = 0.1$ for the step load.

In the case of the step load, the asymptotic solution (4) for the displacement $\varphi_{n,0}(t)$ is derived for $\lambda = 0$ only. However, from Fig. 2a and Fig. 3a it can be seen that it agrees well with the numerical solutions obtained both for $\lambda = 0$ and for $\lambda = 0.1$. This suggests that, in the case of the step load, the asymptotic solution for the displacement $\varphi_{n,0}(t)$ does not depend on the viscosity parameter $\lambda$ and is described by (4). At the same time, the influence of the viscosity parameter $\lambda$ on the behavior of velocity and acceleration is significant.

In Figs. 4a, 4b, numerical and asymptotic solutions (5), (6), (12), and (13) are shown for the step load in the case $\lambda = 0.02$ at the point with coordinates $n = m = 25$.

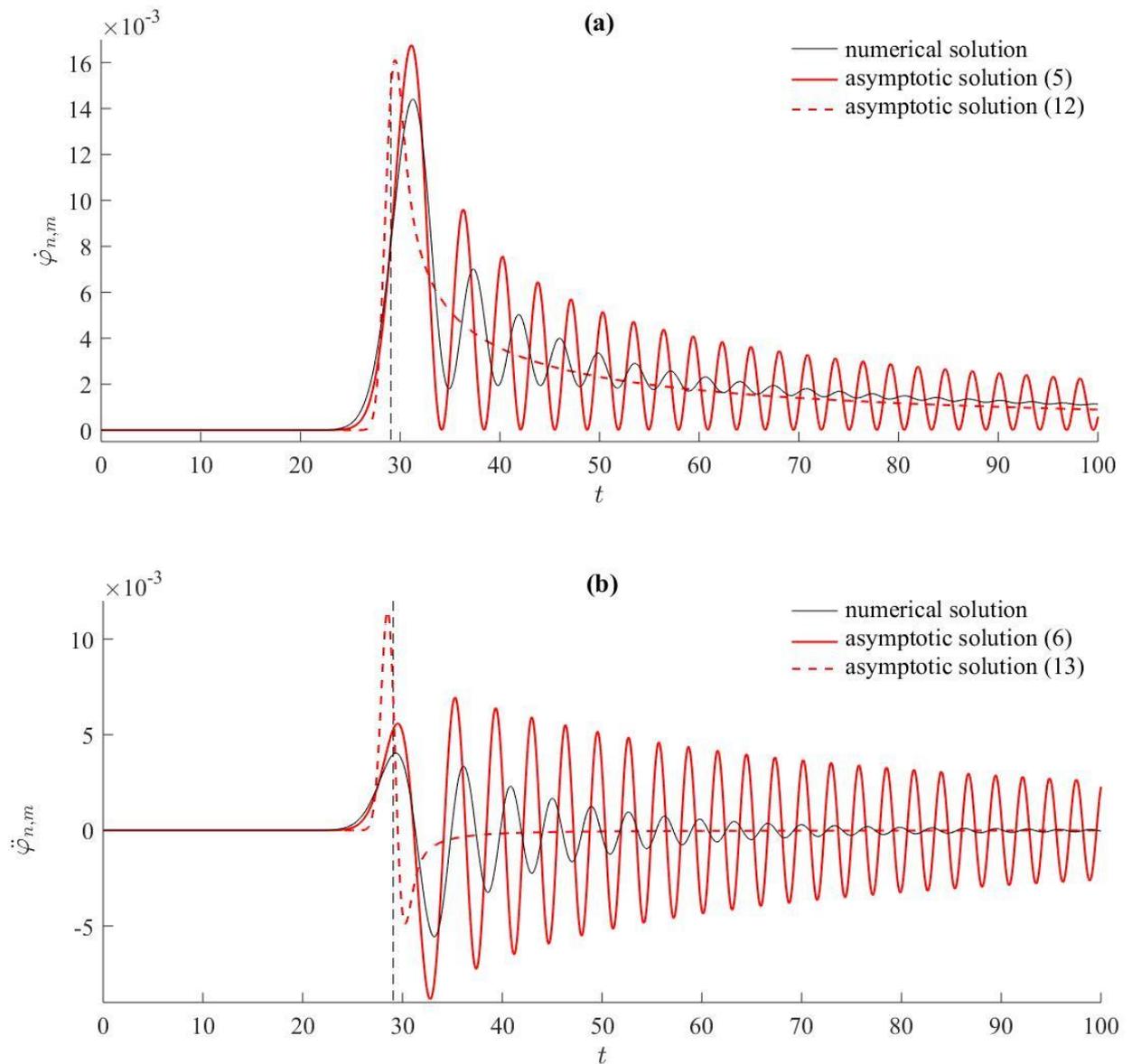

**Fig. 4.** Time dependences of velocities $\dot{\varphi}_{n,m}$, and accelerations $\ddot{\varphi}_{n,m}$ of the mass with coordinates $n = m = 25$ in the case $\lambda = 0.02$ for the step load.

Numerical calculations performed for various values of $\lambda$ showed that, as $\lambda$ increases, the amplitudes of high-frequency oscillations $\dot{\varphi}_{n,m}$ and $\ddot{\varphi}_{n,m}$ decrease and, for $\lambda$ large enough, disappear completely (see Figs. 2 – 4). For example, in Figs. 3b and 3c, it can be seen that for $\lambda = 0.1$ there are no high-frequency oscillations behind the quasifront $n = c_1 t$. Comparison of finite-difference and asymptotic solutions shows that, with increasing $\lambda$, the rate of convergence of the asymptotic solution (12) for velocity $\dot{\varphi}_{n,m}$ to the finite difference solution is significantly higher than the rate of convergence of the asymptotic solution (13) for acceleration $\ddot{\varphi}_{n,m}$ to the finite difference solution. This can be seen in Figs. 3b, 3c and Figs. 4a, 4b. Despite the fact that the asymptotic solutions (5), (6) were obtained for $\lambda = 0$, they can be also used to estimate velocities and accelerations for small values of $\lambda$. Figures 2b, 3b, and 4a show that, in order to estimate the displacement velocity, one can use the asymptotic solution (5) for $0 \le \lambda \le 0.02$ and the asymptotic solution (12) for $\lambda \ge 0.02$. Figures 2c and 3c show that, in order to estimate the accelerations, one can use the asymptotic solution (6) for $\lambda = 0$ and the asymptotic solution (13) for $\lambda \ge 0.1$. Similar figures, not included in this article, show that, in fact, the asymptotic solution (6) can be used to estimate accelerations for $0 \le \lambda \le 0.002$. In the specified intervals for $\lambda$, the asymptotic solutions (4) – (6), (12), and (13) agree well with the numerical solutions both qualitatively and quantitatively.

### 2.2. Gaussian pulse load

In many seismic problems, the study of perturbations caused by an pulse load is of the greatest interest. We study the response of the lattice to the concentrated action of the Gaussian pulse:

$$Q(t) = \exp\left[-(t-4\sigma)^2/(2\sigma^2)\right], \quad \sigma > 0. \tag{17}$$

The Laplace–Fourier transform of the load (17) is as follows:

$$Q^L(p) = \frac{\pi^{1/2}\sigma e^{(\sigma p)^2/2 - 4\sigma p}}{2^{1/2}}[\mathrm{erf}(\frac{4-\sigma p}{2^{1/2}}) + 1], \qquad (18)$$

where $\mathrm{erf}\,x$ is the Gauss error function.

We want to obtain asymptotic solutions for $t \to \infty$ in the following cases: (a) a short Gaussian pulse and $\lambda = 0$; (b) a long Gaussian pulse and arbitrary $\lambda$; (c) a Gaussian pulse of any duration and large values of $\lambda$. Values of $\sigma$ for which the Gaussian pulse is considered short or long, as well as values of $\lambda$ which are considered large, will be specified below. Note that in case (a), both high-frequency oscillations and low-frequency oscillations are excited in the lattice; in case (b), only low-frequency oscillations are excited in the lattice; in case (c), even if high-frequency oscillations are excited, the viscous dampers gradually dampen them, and over time only low-frequency perturbations remain in the lattice. Therefore, in case (a), we use the method proposed in [12]. It allows us to take into account both low-frequency and high-frequency oscillations. In cases (b) and (c), we are looking for only the low-frequency asymptotic solution, using Slepyan's method [17] for this.

Let us start finding an asymptotic solution for case (a). We take the duration of the Gaussian pulse (17) equal to $8\sigma$, and we call the pulse itself short if $\sigma < 1/(8c_1)$.

We find the inverse Laplace and Fourier transforms of the function $\varphi^{LF_n}$ given by (3). To do this, we use the method proposed in [12]. As a result, from (3) and (18) we get the following asymptotic formula for displacements $\varphi_{n,0}$ as $t \to \infty$ and $n \to \infty$:

$$\varphi_{n,0}(t) \sim \frac{\pi^{1/2}\sigma}{2^{1/2}c_1} J_n^2[c_1(t-4\sigma)]H(t-4\sigma). \qquad (19)$$

Differentiating (19) with respect to time, we find the asymptotics of the velocity $\dot{\varphi}_{n,0}$ and acceleration $\ddot{\varphi}_{n,0}$ of the lattice masses as $t \to \infty$ and $n \to \infty$:

$$\dot{\varphi}_{n,0}(t) \sim (2\pi)^{1/2}\sigma J_n[c_1(t-4\sigma)]J_n'[c_1(t-4\sigma)]H(t-4\sigma), \qquad (20)$$

$$\ddot{\varphi}_{n,0}(t) \sim (2\pi)^{1/2}\sigma c_1\{J'_n[c_1(t-4\sigma)]J'_n[c_1(t-4\sigma)]+J_n[c_1(t-4\sigma)]J''_n[c_1(t-4\sigma)]\}H(t-4\sigma). \qquad (21)$$

Asymptotic solutions (19) – (21) take into account both long-wave perturbations ($q_x$, $q_y \to 0$) and short-wave perturbations ($q_x = 0$, $q_y = \pi$ or $q_x = \pi$, $q_y = 0$).

Using (7) and (8), we rewrite (19) – (21) in a different form:

$$\varphi_{n,0}(t) \sim \frac{\pi^{1/2}\sigma[\mathrm{Ai}(\kappa)]^2}{2^{1/2}c_1[c_1(t-4\sigma)/2]^{2/3}}H(t-4\sigma), \qquad (22)$$

$$\dot{\varphi}_{n,0}(t) \sim -\frac{2^{3/2}\pi^{1/2}\sigma\mathrm{Ai}(\kappa)\mathrm{Ai}'(\kappa)}{c_1(t-4\sigma)}H(t-4\sigma), \qquad (23)$$

$$\ddot{\varphi}_{n,0}(t) \sim \frac{(2\pi)^{1/2}\sigma c_1\{[\mathrm{Ai}'(\kappa)]^2 + \mathrm{Ai}(\kappa)\mathrm{Ai}''(\kappa)\}}{[c_1(t-4\sigma)/2]^{4/3}}H(t-4\sigma), \qquad (24)$$

where

$$\kappa = \frac{n-c_1(t-4\sigma)}{[c_1(t-4\sigma)/2]^{1/3}} \quad \text{and} \quad \sigma > 0. \qquad (25)$$

When deriving (24) from (21), we used the formula

$$J''_n(c_1 t) \sim \frac{\mathrm{Ai}''(\kappa)}{(c_1 t/2)},$$

which is obtained by differentiating (8) and passing to the limit as $t \to \infty$ and $n \to \infty$.

Formulas (22) – (25) show that, in the case when there is no viscosity ($\lambda = 0$) and the loading is carried out by a short Gaussian pulse, the following effects are observed in the lattice as $t \to \infty$ ($n \to \infty$): the maximum amplitude of the displacements of the masses of the lattice decreases as $t^{-2/3}$ ($n^{-2/3}$); the maximum amplitude of the velocities of the masses of the lattice decreases as $t^{-1}$ ($n^{-1}$); the maximum amplitude of the accelerations of the masses of the lattice decreases as $t^{-4/3}$ ($n^{-4/3}$); and the quasi-front zone expands as $t^{1/3}$ ($n^{1/3}$).

In Figs. 5a – 5c, numerical and asymptotic solutions (19) – (21) are shown for the short Gaussian pulse load ($\sigma = 0.1$) in the case $\lambda = 0$ at the point with coordinates $n = m = 25$. For this, asymptotic solutions (19) – (21) are rewritten in terms of the radial coordinate $r = \sqrt{n^2 + m^2}$. In Fig. 5, the vertical dashed lines correspond to the arrival time of the quasi-front of the longitudinal wave, $t = r/c_1 + 4\sigma$.

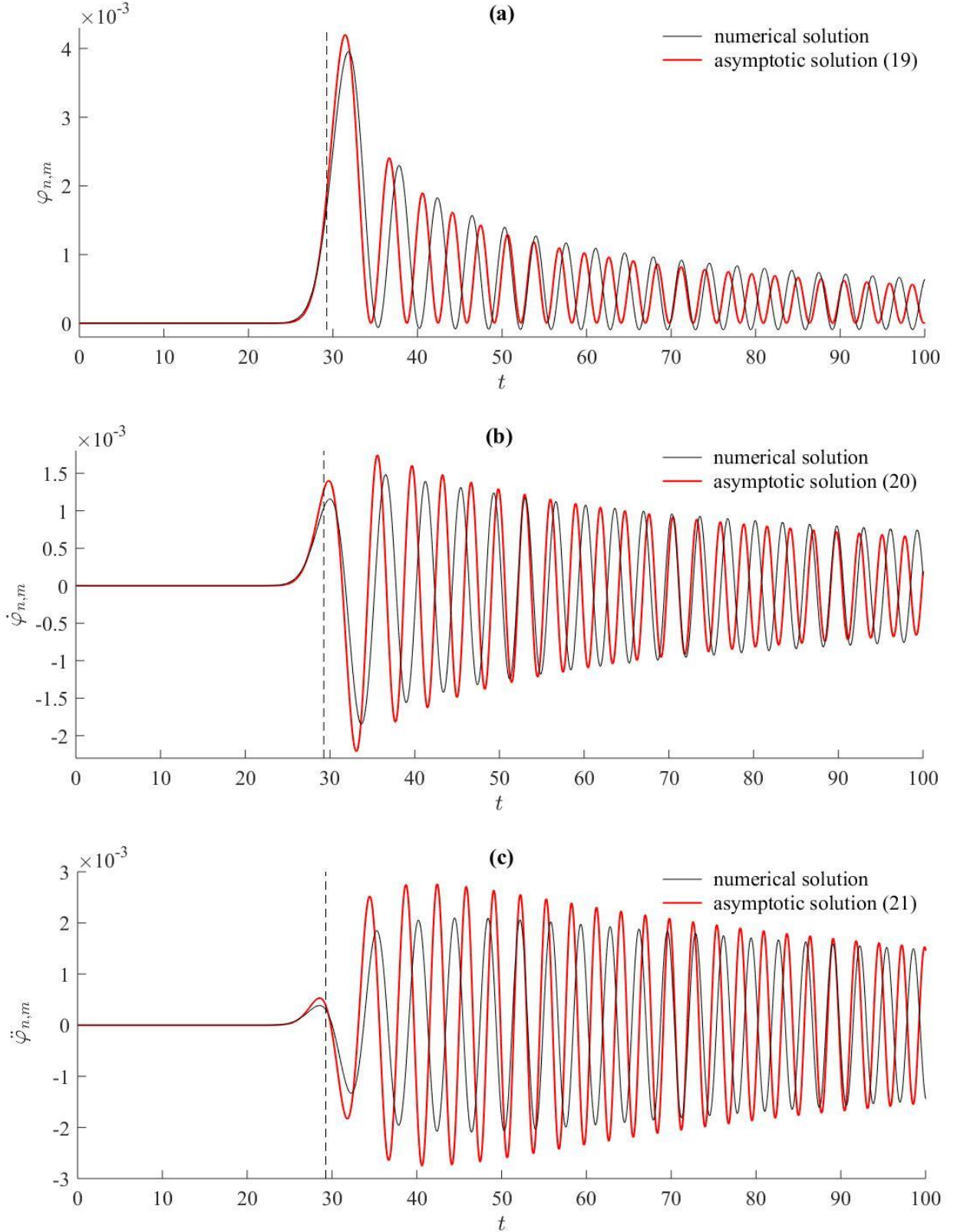

**Fig. 5.** Time dependences of displacements $\varphi_{n,m}$, velocities $\dot{\varphi}_{n,m}$, and accelerations $\ddot{\varphi}_{n,m}$ of the mass with coordinates $n = m = 25$ in the case $\lambda = 0$ for the Gaussian pulse load ($\sigma = 0.1$).

Fig. 5 shows that in the case of the short Gaussian pulse load in the lattice without viscous dampers ($\lambda = 0$) there are high-frequency oscillations behind the quasi-front $n = c_1(t - 4\sigma)$. We saw the same effect in Fig. 2b and Fig. 2c for a step load in case of $\lambda = 0$.

Numerical calculations carried out for $\lambda = 0$ and various values of $\sigma$ allow us to conclude that the correspondence between the maximum amplitudes of the asymptotic solutions (19) – (21) and finite-difference solutions is acceptable for $\sigma \leq 0.1$, but worsens with the growth of $\sigma$.

Let's start finding a low-frequency asymptotic solution for cases (b) and (c), when the Gaussian pulse is long and $\lambda$ is arbitrary or the duration of the Gaussian pulse is arbitrary and the value of $\lambda$ is large. Which $\lambda$ are considered large and for which $\sigma$ the Gaussian pulse is considered long, we will determine later from comparison of asymptotic and finite difference solutions.

To find the inverse Fourier transform with respect to $q_x$ and the inverse Laplace transform with respect to $p$ of the function $\dot{\varphi}_0^{LF_n}$ in (3), we use Slepyan's asymptotic method [17] of joint inversion of these transforms in a vicinity of the ray $n = c_1 t$. Taking into account (18), we obtain the following asymptotic solutions for the displacements $\varphi_{n,0}$ of the lattice masses, their velocities $\dot{\varphi}_{n,0}$, and accelerations $\ddot{\varphi}_{n,0}$:

$$\varphi_{n,0}(t) = \frac{2^{1/4} \sigma \Phi_1(\kappa)}{3\pi(t - 4\sigma)^{1/2}[\lambda(t - 4\sigma) + \sigma^2]^{1/4}} H(t - 4\sigma), \quad (26)$$

$$\dot{\varphi}_{n,0}(t) = \frac{2^{3/4} \sigma \Phi_2(\kappa)}{3\pi(t - 4\sigma)^{1/2}[\lambda(t - 4\sigma) + \sigma^2]^{3/4}} H(t - 4\sigma), \quad (27)$$

$$\ddot{\varphi}_{n,0}(t) = \frac{2^{5/4} \sigma \Phi_3(\kappa)}{3\pi(t - 4\sigma)^{1/2}[\lambda(t - 4\sigma) + \sigma^2]^{5/4}} H(t - 4\sigma). \quad (28)$$

Here

$$\kappa = \frac{n - c_1(t - 4\sigma)}{(3/4)^{1/2}[\lambda(t - 4\sigma) + \sigma^2]^{1/2}}, \tag{29}$$

the functions $\Phi_1(\kappa)$ and $\Phi_2(\kappa)$ are defined by (14) and (15), and the function $\Phi_3(\kappa)$ is defined by the formula

$$\begin{aligned}\Phi_3(\kappa) &= \int_0^\infty e^{-z^2} z^{3/2} \sin(z\kappa - \pi/4) dz = \\ &= 2^{-7} \pi e^{-\eta} |\kappa|^{5/2} \left\{ 4I_{-1/4}(\eta) - 5I_{3/4}(\eta) + I_{-5/4}(\eta) - \operatorname{sgn}(\kappa)\left[4I_{1/4}(\eta) - 5I_{-3/4}(\eta) + I_{5/4}(\eta)\right]\right\}.\end{aligned} \tag{30}$$

where $\eta = \kappa^2/8$.

Using formulas $I_{1/4}(\varsigma) - 2\varsigma\left[I_{-3/4}(\varsigma) - I_{5/4}(\varsigma)\right] = 0$ and $I_{-1/4}(\varsigma) - 2\varsigma\left[I_{3/4}(\varsigma) - I_{-5/4}(\varsigma)\right] = 0$ which are valid for all $\varsigma$, we rewrite (30) in the form

$$\Phi_3(\kappa) = 2^{-5} \pi e^{-\eta} |\kappa|^{1/2} \left\{ (\kappa^2 - 2)I_{-1/4}(\eta) - \kappa^2 I_{3/4}(\eta) - \operatorname{sgn}(\kappa)\left[(\kappa^2 - 2)I_{1/4}(\eta) - \kappa^2 I_{-3/4}(\eta)\right]\right\}.$$

Comparing the last formula with (14), (15), we obtain $\Phi_3(\kappa) \sim -\Phi_1(\kappa)/4 + \kappa\Phi_2(\kappa)/2$ or

$$\ddot{\varphi}_{n,0}(t) \sim \frac{2^{3/2} 3^{-1/2} \dot{\varphi}_{n,0}(t)[n - c_1(t - 4\sigma)] - \varphi_{n,0}(t)}{2[\lambda(t - 4\sigma) + \sigma^2]}.$$

Thus, in cases (b) and (c) the asymptotic solution for acceleration is expressed in terms of the asymptotic solutions for displacement and velocity.

Formulas (26) – (29) show that, in cases (b) and (c), the following effects are observed as $t \to \infty$ ($n \to \infty$): the maximum amplitude of the displacements $\varphi_{n,0}$ of the masses of the lattice decreases as $t^{-3/4}$ ($n^{-3/4}$); the maximum amplitude of their velocities $\dot{\varphi}_{n,0}$ decreases as $t^{-5/4}$ ($n^{-5/4}$); the maximum amplitude of their accelerations $\ddot{\varphi}_{n,0}$ decreases as $t^{-7/4}$ ($n^{-7/4}$); and the quasi-front zone expands as $t^{1/2}$ ($n^{1/2}$). In addition, Eqs. (26) – (29) show that the maximum amplitude of the displacements $\varphi_{n,0}$ of the masses of the lattice decreases as $\lambda^{-1/4}$ with an increase in the viscosity parameter $\lambda$; the maximum amplitude of their velocities $\dot{\varphi}_{n,0}$ decreases as $\lambda^{-3/4}$; the maximum amplitude of their accelerations $\ddot{\varphi}_{n,0}$ decreases as $\lambda^{-5/4}$; and the quasi-front zone expands as $\lambda^{1/2}$.

In Figs. 6a – 6c, numerical and asymptotic solutions (26) – (28) are shown for the Gaussian pulse load ($\sigma = 5$) in the case $\lambda = 0.1$ at the point with coordinates $n = m = 25$. For this, asymptotic solutions (26) – (28) are rewritten in terms of the radial coordinate $r = \sqrt{n^2 + m^2}$. In Fig. 6, the vertical dashed lines correspond to the arrival time of the quasi-front of the longitudinal wave, $t = r/c_1 + 4\sigma$. Note that this time depends on $\sigma$ and, for large values of $\sigma$, can significantly exceed the value $r/c_1$. The comparison of Fig. 5 and Fig. 6 confirms this conclusion.

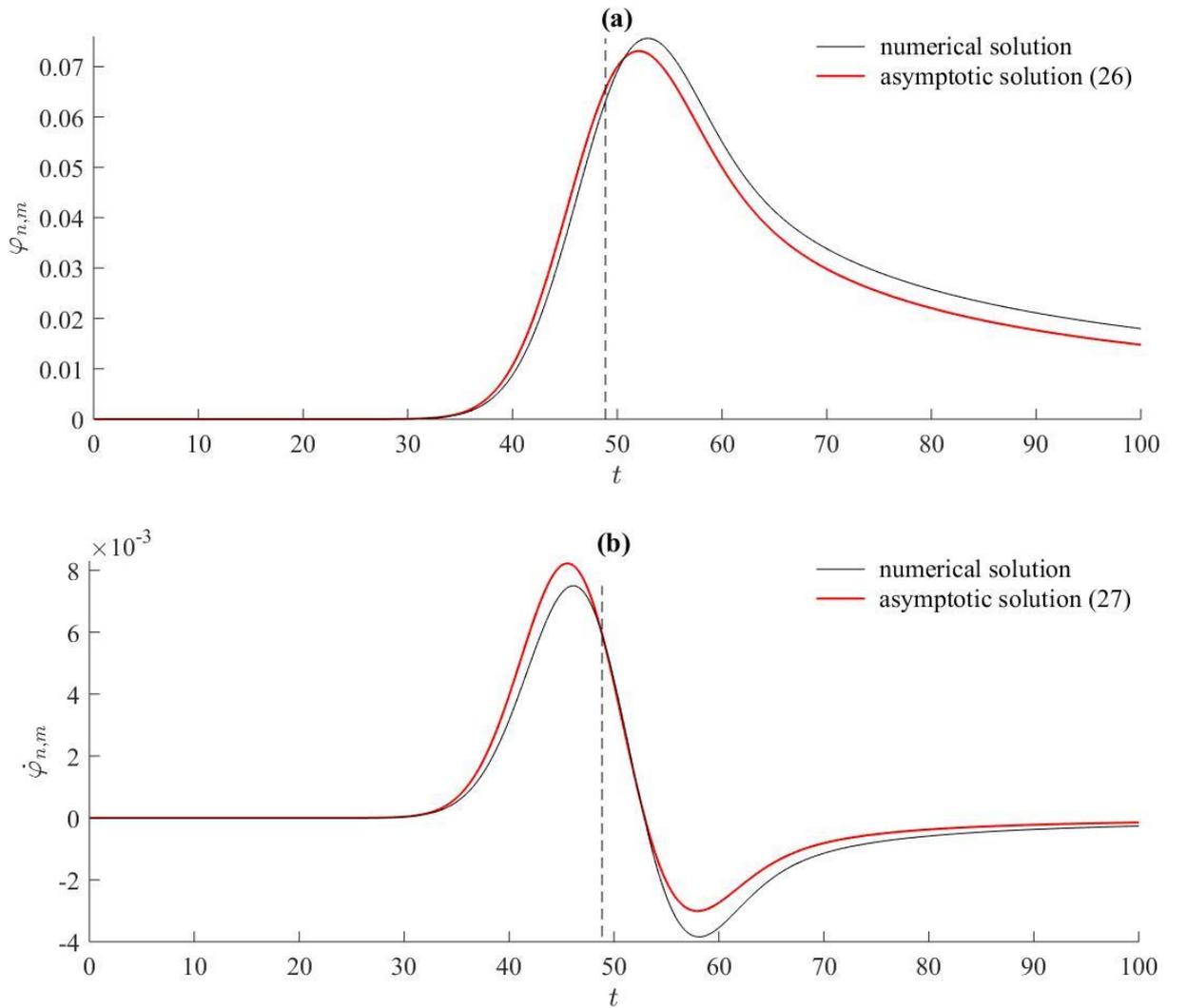

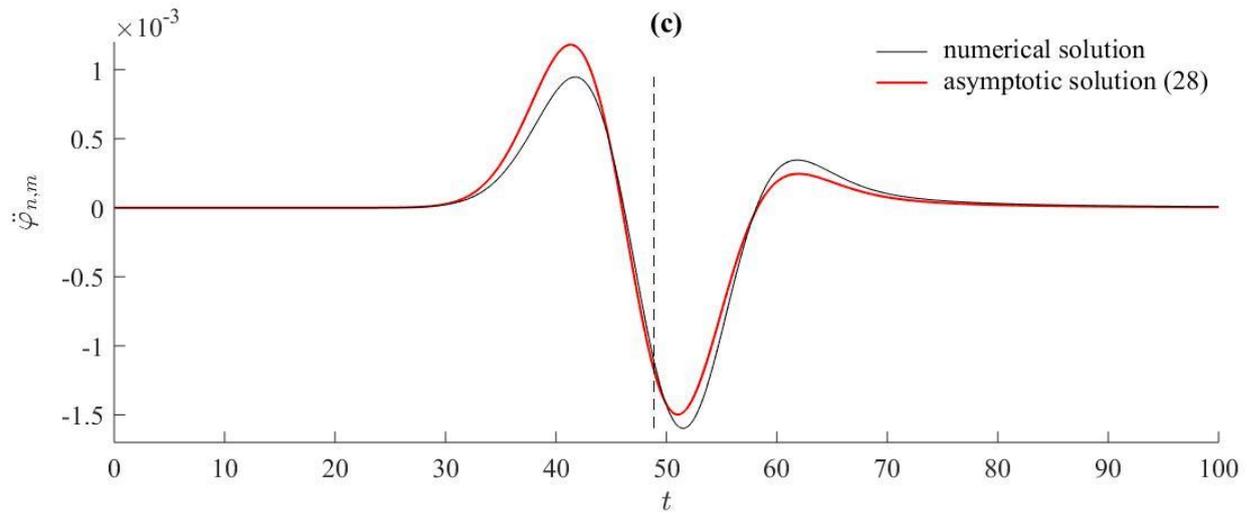

**Fig. 6.** Time dependences of displacements $\varphi_{n,m}$, velocities $\dot{\varphi}_{n,m}$, and accelerations $\ddot{\varphi}_{n,m}$ of the mass with coordinates $n=m=25$ in the case $\lambda=0.1$ for the Gaussian pulse load ($\sigma=5$).

Let us determine the values of $\sigma$ and $\lambda$ for which it is possible to use the high-frequency asymptotic solution (19) – (21), and for which the low-frequency asymptotic solution (26) – (28). That is, we determine which values of $\sigma$ and $\lambda$ correspond to each of the cases (a), (b), (c). We do this below by comparing the mentioned asymptotic solutions with finite difference solutions.

Despite the fact that in case (a) the asymptotic solutions (19)–(21) were obtained for $\lambda=0$ only, numerical calculations show that these solutions can be also used for small values of $\lambda$. Namely, if $\sigma<1/(8c_1)$, then to describe the displacements, one can use (19) provided that $\lambda \leq 0.01$; to describe the velocities, one can use (20) provided that $\lambda \leq 0.05$; and to describe the accelerations, one can use (21) provided that $\lambda \leq 0.001$. We consider these values of $\lambda$ to be small.

In case (b), if $\lambda=0$, then the asymptotic solution (26) can be used to describe the displacements of the masses $\varphi_{n,m}$ provided that $\sigma \geq 10/(8c_1)$ and the asymptotic solutions (27), (28) can be used to describe their velocities $\dot{\varphi}_{n,m}$ and accelerations $\ddot{\varphi}_{n,m}$ provided that $\sigma \geq 30/(8c_1)$. If, $\lambda=0.05$ then the asymptotic solution (26) can be used to describe the displacements of the masses $\varphi_{n,m}$ provided that $\sigma \geq 2/(8c_1)$ and the asymptotic solutions (27), (28) can be used to describe their velocities $\dot{\varphi}_{n,m}$ and accelerations $\ddot{\varphi}_{n,m}$ provided that $\sigma \geq 30/(8c_1)$. Under these restrictions, we consider the Gaussian pulse to be long. Note that these restrictions on $\sigma$ turn out to be different for different values of $\lambda$ and for different functions $\varphi_{n,m}$, $\dot{\varphi}_{n,m}$, $\ddot{\varphi}_{n,m}$.

In case (c), if $\sigma \leq 1/(8c_1)$, then the asymptotic solution (26) can be used to describe the displacements $\varphi_{n,m}$ provided that $\lambda \geq 0.05$ and the asymptotic solutions (27), (28) can be used to describe the velocities $\dot{\varphi}_{n,m}$ and accelerations $\ddot{\varphi}_{n,m}$ provided that $\lambda \geq 0.2$. We consider the above values of $\lambda$ to be large for the corresponding function $\varphi_{n,m}$, $\dot{\varphi}_{n,m}$, or $\ddot{\varphi}_{n,m}$. Finally, we note that for $\sigma > 1/(8c_1)$ the above inequalities on $\lambda$ also allow one to use formulas (26) – (28), though in fact they can be replaced by slightly less restrictive ones.

**Conclusion**

The propagation of transient longitudinal perturbations in 2D square lattices of point masses connected by Voigt elements in axial and diagonal directions is studied. Concentrated transient loads orthogonal to the lattice plane are considered.

Numerical and asymptotic solutions are obtained for a long time from the beginning of the process or at a large distance from the place of the application of the load for two types of the time dependence of the amplitude of the load, namely, for the Heaviside step function and for the Gaussian pulse.

The viscosity of the dampers and the loading parameters are determined for which a high-frequency or low-frequency wave process is formed in the lattice.

Using asymptotic methods, the degree of the attenuation of the amplitudes of perturbations in the lattice with increasing time or distance is determined. The dependence of the degree of attenuation on the viscosity parameter and the type of applied load is studied.

Asymptotic and finite-difference solutions are compared with each other for lattices subjected to various loads and having various viscosities of the dampers. It is shown that they agree with each other qualitatively and quantitatively. This agreement occurs at a finite exposure time or, which is the same, at a finite distance from the place of application of the load. The limits of applicability of asymptotic solutions are determined depending on the parameters of the problem.

**Declaration of Competing Interest**

The author declares that she has no known competing financial interests or personal relationships that could have appeared to influence the work reported in this paper.

**Acknowledgement**